\documentclass[11pt,a4paper]{article}

\usepackage{graphics,epsfig,theorem,latexsym,amssymb,amsmath, amsfonts}
\usepackage{times,epic,eepic}

\usepackage{pifont}

\usepackage{cite}
\usepackage{url}

\usepackage{color}

\usepackage{newlfont}

 \newenvironment{keywords}{ \noindent {\small\bf Key Words}:}{ }

\def\G1{\hbox{$\displaystyle{\mbox{\ding{172}}}$}}

\def\bd{\begin{description}}
\def\ed{\end{description}}

\def\beq{\begin{equation}}
\def\eeq{\end{equation}}
\def\bea{\begin{eqnarray}}
\def\eea{\end{eqnarray}}
\def\beas{\begin{eqnarray*}}
\def\eeas{\end{eqnarray*}}

\theoremstyle{remark}

\begin{document}
%\begin{article}

\title{ The Olympic Medals Ranks, lexicographic ordering \\ and numerical infinities  }

%% Author name
%\newcommand{\nms}{\normalsize}
\newcommand{\nms}{\normalsize}
\author{\\ {   \bf Yaroslav D. Sergeyev$^{1-3}$
    }\\[-4pt]
    \\
      \nms $^1$\small{University of Calabria, Rende (CS), Italy}\\ [-4pt]
        \nms $^2$\small{N.I.~Lobatchevsky State University,
  Nizhni Novgorod, Russia}\\[-4pt]
       \nms $^3$\small{Institute of High Performance Computing
   and Networking}\\[-4pt]
       \nms  \small{of the National Research Council of Italy, Rende (CS), Italy}\\[-4pt]
         \nms {\tt \small{ yaro@si.dimes.unical.it }}
}

%\titlerunning{Higher order numerical differentiation  on the Infinity Computer}
%\authorrunning{Yaroslav D. Sergeyev}

\date{}

\maketitle

\vspace*{-5mm}
 \begin{abstract}

Several ways used to rank countries with respect to medals won
during Olympic Games are discussed. In particular, it is shown that
the unofficial rank used by the Olympic Committee   is the only rank
that does not allow one to use a numerical counter for ranking --
this rank uses the lexicographic ordering to rank countries: one
gold medal is more precious than any number of silver medals and one
silver medal is more precious than any number of bronze medals. How
can we quantify what do these words, \textit{more precious}, mean?
Can we introduce a counter that for \textit{any} possible number of
medals would allow us to compute a \textit{numerical} rank of a
country using the number of gold, silver, and bronze medals in such
a way that the higher resulting number would put the country in the
higher position in the rank? Here we show that it is impossible to
solve this problem using
  the positional numeral system with \textit{any}
finite base. Then we demonstrate that this  problem can be easily
solved by applying numerical computations with recently developed
actual infinite numbers. These computations can be done on a new
kind of a computer -- the recently patented Infinity Computer. Its
working software prototype is described briefly and examples of
 computations are given. It is shown that the new way of counting can
be used in all situations where the lexicographic ordering is
required.
 \end{abstract}

%\vspace{5mm}


%\subclass{65K05, 90C26, 90C56}

%\vspace{5mm}

\section*{Olympic  medals ranks}
\label{s_m1}

The International Olympic Committee (IOC) does not produce any
official   ranking of  countries participating at the Olympic Games.
However,   the IOC  publishes tables showing  medals won by athletes
representing each country participating at the Games.   The
convention used by the IOC   to order the countries in this
unofficial rank is the following. First, countries are sorted by the
number of gold medals won. If the number of gold medals is the same,
the number of silver medals is taken into consideration, and then
the number of bronze ones. If two countries have an equal number of
gold, silver, and bronze medals, then equal ranking is given and the
countries are listed alphabetically by their IOC country code (for
instance, in the 2010 Winter Olympics held in Vancouver,  China and
Sweden have won   5 gold,  2 silver, and  4 bronze medals; both
countries have the 7$^{th}$ place in the rank but China is higher in
the table). Table 1 shows countries sorted by using this   rank at
Sochi 2014 Olympic Games (the first ten countries). This rank will
be called R1 hereinafter.

\begin{table}[t]
\caption{The International Olympic Committee  unofficial medal rank
at Sochi 2014   (the first ten countries).}
\begin{center} \scriptsize \label{table1}
\begin{tabular}{@{\extracolsep{\fill}}|c|c|c|c|c|c|}\hline
  &   &   &   &      \vspace{-2mm}\\
Rank R1 & Country & Gold & Silver &  Bronze  \\
  &   &   &   &      \vspace{-2mm} \\
 \hline
1 & Russian Federation&  13 & 11  & 9     \\
2 & Norway &  11 & 5  & 10     \\
3 &  Canada  &  10 & 10  & 5     \\
4 &  United States &  9 & 7  & 12     \\
5 & Netherlands  &  8 & 7  & 9     \\
6 & Germany   &  8 & 6  & 5     \\
7 & Switzerland  &  6 & 3  & 2     \\
8 & Belarus   &  5 & 0  & 1    \\
9 & Austria &  4 & 8  & 5    \\
10 & France &  4 & 4  & 7     \\
\hline
\end{tabular}
\end{center}
\end{table}

However, there exist several methods to rank countries (some of them
are illustrated in Tables~2  and~3 showing best 10 countries for
each rank, for more countries see, e.g., \cite{Florida}). First, in
many countries  the total number of Olympic medals earned by
athletes representing each country is very popular. This rank
(called hereinafter R2) gives equal ratings to gold, silver, and
bronze medals. So, if a country $A$ has won $g_A$ gold, $s_A$
silver, and $b_A$ bronze medals than its rank is just
\[
R2(A) = g_A + s_A + b_A.
\]
Since this rank gives the same weight to gold, silver, and bronze
medals, there are several proposals to improve this way of counting
by introducing weights for medals. For instance, the Fibonacci
weighted point system (this method is shown in Table~2 as R3) uses
the following weights: gold costs 3 points, silver 2 points, and
bronze 1 point, these weights are called 3:2:1 system. Thus, it
follows
\[
R3(A) = 3g_A + 2s_A +  b_A.
\]
It can be seen from Table 2 that Norway and United States have the
same rank R3 but in the rank Norway has a higher position because it
has won more golds (the same situation holds for Switzerland  and
Sweden). In order to make gold medals more precious, the exponential
weighted point system assigns 4 points to gold, 2 points to silver,
and 1 point to bronze -- 4:2:1 system. The variation employed by the
British press during the Olympic Games in London in 1908 used the
weights 5:3:1. There exist also systems 5:3:2, 6:2:1, 10:5:1, etc.

\begin{table}[t]
\caption{Medal ranks counting total number of won medals per country
and weighted total sum (system 3:2:1).}
\begin{center} \scriptsize \label{table2}
\begin{tabular}{@{\extracolsep{\fill}}|c|cc|cc| }
\hline
  & \multicolumn{2}{c|}{  } &
\multicolumn{2}{c|}{  }      \vspace{-2mm}\\
 $N$   & \multicolumn{2}{c|}{Total medals (R2) } &
\multicolumn{2}{c|}{ Weighted total medals (R3)}  \\
 & \multicolumn{2}{c|}{  } &
\multicolumn{2}{c|}{  }      \vspace{-2mm}\\
\hline
1 & Russian Federation  & 33 & Russian Federation  & 70  \\
2 & United States       & 28 & Canada              & 55  \\
3 & Norway              & 26 & Norway              & 53   \\
4 & Canada              & 25 & United States       & 53  \\
5 & Netherlands         & 24 & Netherlands         & 47   \\
6 & Germany             & 19 & Germany             & 41  \\
7 & Austria             & 17 & Austria             & 33     \\
8 & France              & 15 & France              & 27   \\
9 & Sweden              & 15 & Switzerland         & 26   \\
10 & Switzerland        & 11 & Sweden              & 26  \\
\hline
\end{tabular}
\end{center}
\end{table}

\begin{table}[t]
\caption{Medal ranks counting  total medals per 10 million people
and total medals  per \$100 billion of the gross domestic product.}
\begin{center} \scriptsize \label{table3}
\begin{tabular}{@{\extracolsep{\fill}}|c|cr|cc| }
\hline
  & \multicolumn{2}{c|}{  } &
\multicolumn{2}{c|}{  }      \vspace{-2mm}\\
 $N$   & \multicolumn{2}{c|}{Total   medals per $10^7$ people (R4) } &
\multicolumn{2}{c|}{ Total medals per \$100 billion of  GDP (R5)}  \\
 & \multicolumn{2}{c|}{  } &
\multicolumn{2}{c|}{  }      \vspace{-2mm}\\
\hline
1 & Norway         & 51.8 & Slovenia            & 17.7  \\
2 & Slovenia       & 38.9 & Latvia              & 14.1  \\
3 & Austria        & 20.1 & Belarus             & 9.5   \\
4 & Latvia         & 19.8 & Norway              & 5.2  \\
5 & Sweden         & 15.8 & Austria             & 4.3   \\
6 & Netherlands    & 14.3 & Czech Republic      & 4.1  \\
7 & Switzerland    & 13.8 & Netherlands         & 3.1    \\
8 & Finland        &  9.2 & Sweden              & 2.9   \\
9 & Czech Republic & 7.6  & Finland             & 2.0   \\
10 & Canada        &  7.2 & Switzerland         & 1.7  \\
\hline
\end{tabular}
\end{center}
\end{table}

  Ranks using completely different ideas should be also mentioned. For
instante, it is possible to count   all the   medals won (weighted
or not) counting separately the medals for each individual athlete
in team sports. Another idea consists of the usage of an improvement
rank that is based on percentage improvement reached by countries
with respect to the previous Games results. There exist ranks built
in comparison to expectations predictions. Among them, there can be
mentioned predictions using previous results (during the Games or
other competitions) and predictions using economics, population, and
a range of other.

Another interesting proposal is to calculate the rank dividing the
number of medals   by the population of the country. The column R4
in Table~3 shows total medals won by a country divided by 10 million
people. While criteria R1 -- R3 give similar results, criterion R4
puts different countries, mainly having relatively small
populations, on the top of the rank. In fact, Norway, that has won
26 medals and has a population of approximately 5 million people is
the best in this rank. In general, the countries that top the list
have also small populations in comparison, for instante, with United
States and Russian Federation. The number of medals per \$100
billion of the gross domestic product (GDP) of the country (this
rank is called R5 in Table~2)   also favors smaller countries.

In this paper, we do not discuss advantages and disadvantages of
various ranks. We consider a purely mathematical problem regarding
the following difference that distinguishes the   unofficial
International Olympic Committee rank R1 and the other ranks R2 --
R5. In fact, while ranks R2 -- R5 produce numerical coefficients for
each country that allow one to order the countries, rank R1 does not
produce any resulting number that can be used for this purpose. This
rank uses the so called \textit{lexicographic ordering} called so
because it is used in dictionaries to order words: first words are
ordered with respect to the first symbol in the word, then with
respect to the second one, and so on. In working with the rank R1 we
have words that consists of three symbols $g_A, s_A, b_A$ and,
therefore, their length $w=3$.

In this paper, we show that, as it happens for ranks R2 -- R5,  it
is possible to propose a procedure for computing rank R1 numerically
for each country and for any number of medals. Moreover, it is shown
that the introduced way of computation can be generalized from words
consisting of three symbols to words having a general finite length
$w$  and used in situations that require the lexicographic ordering.

\section*{How to compute the rank R1 for any   number of medals?}
\label{s_m2}

  Evidently, in the rank R1 gold medals are more precious than
silver ones that in their turn are better than the bronze ones. An
interesting issue in this way of counting consists of the following
fact. Let us consider Belarus   and Austria that occupy the $8^{th}$
and $9^{th}$ positions, respectively. Belarus has 5 gold medals and
Austria only 4. In spite of the fact that Austria has 8 silver
medals and Belarus none of them, this fact is not taken into
consideration. Austria could have \textit{any} number of silver
medals but one gold medal of Belarus will be more important than all
these silver medals.

Can we quantify what do these words, \textit{more important}, mean?
Can we introduce in a way a counter that would allow us to compute a
numerical rank of a country using the number of gold, silver, and
bronze medals in such a way that the higher resulting number would
put the country in the higher position in the rank? Moreover, we
wish to construct a numerical counter that would work for
\textit{any} number of medals. This counter should work also in
situations when the number of medals that can be won is not known a
priori.

More formally, we would like to introduce a number $n(g_A,s_A, b_A)$
where $g_A$ is the number of gold medals,   $s_A$ is the number of
silver medals, and   $b_A$ is the number of the bronze ones won by a
country $A$.    This number should be calculated in such a way that
for countries $A$ and $B$ it follows that
 \beq
n(g_A,s_A, b_A) > n(g_B,s_B, b_B), \,\,  \mbox{if} \left\{
\begin{array}{l} g_A
>
g_B,\\
 g_A =  g_B,
s_A  > s_B,\\
g_A =  g_B, s_A  = s_B, b_A  > b_{0}.
\end{array} \right.
 \label{medals1}
       \eeq
In addition,   $n(g_A,s_A, b_A)$ should be introduced under
condition that the number $K > \max \{ g_A,s_A, b_A \}$ being an
upper bound for the number of medals of each type that can be won by
each country is unknown.

In order to calculate $n(g_A,s_A, b_A)$, let us try to give weights
to $g_A,\, s_A,$ and  $b_A$ as it is done in the positional numeral
system with a base $\beta$:
 \beq
n(g_A,s_A, b_A)=
 g_A \beta^{2} +   s_A \beta^1 +b_A \beta^0 = g_As_A b_A \label{medals6}
       \eeq
For instance, in the decimal positional numeral system with $\beta =
10$ the record
 \beq
n(g_A,s_A, b_A)=
 g_A 10^{2} +   s_A 10^1 +b_A 10^0 = g_As_A b_A \label{medals2}
       \eeq
gives us the rank of the country $A$. However, it can be seen
immediately that this way of reasoning does not solve our problem,
since it does not satisfy condition (\ref{medals1}). In fact, if a
country will have more than 11 silver medals, then formula
(\ref{medals2}) implies that these medals are more important than
one gold. For instance, the data
 \beq
  g_A=2, s_A=0, b_A=0, \hspace{1cm}
g_B=1,s_B=11, b_B=0.
 \label{medals3}
       \eeq
give us
\[
n(g_A,s_A, b_A)= 2 \cdot 10^{2} +  0 \cdot 10^1 + 0 \cdot 10^0 = 200
<
\]
\[
n(g_B,s_B, b_B)= 1 \cdot 10^{2} +  11 \cdot 10^1 + 0 \cdot 10^0 =
210,
\]
i.e., condition (\ref{medals1}) is not satisfied.

Remind that we wish to construct a numerical counter that would work
for \textit{any} number of medals, in other words, we suppose that
countries can win any number of medals and this number is unknown
for us. Then it is easy to see that situations can occur where the
positional system will not satisfy (\ref{medals1}) not only for the
base $\beta=10$ but also for any finite $\beta$. This can happen if
one of the countries will have more than $\beta$ silver (or bronze)
medals.

Thus, the contribution of 1 gold medal in the computation of
$n(g_A,s_A, b_A)$ should be larger than the contribution of
\textit{any}   number, $s_A$, of silver medals, i.e., it should be
\emph{infinitely larger}. Analogously, the contribution of 1 silver
medal should be infinitely larger than the contribution of any
finite number of bronze medals.

Unfortunately, it is difficult to make numerical computations with
infinity (symbolic ones can be done using non-standard analysis
approach, see \cite{Robinson}) since in the traditional calculus
$\infty$ absorbs any finite quantity and we have, for instance,
 \beq
\infty + 1= \infty,    \hspace{1cm}    \infty + 2 = \infty.
    \label{Riemann30}
 \eeq

\section*{A numerical calculator of the   rank R1
involving infinities}
 \label{s_m3}

In order to   construct a numerical calculator of medal rank
involving infinite numbers, let us    remind the difference
between \emph{numbers} and \emph{numerals}: a \textit{numeral} is
a symbol or group of symbols that represents a \textit{number}.
The difference between them is the same as the difference between
words and the things they refer to. A \textit{number} is a concept
that a \textit{numeral} expresses. The same number can be
represented by different numerals. For example, the symbols `7',
`seven', and `VII' are different numerals, but they all represent
the same number.

Different numeral systems can represent different numbers. For
instance, Roman numeral system is not able to represent zero and
negative numbers. There exist even weaker numeral systems.
Recently (see \cite{Gordon}) a study on a numeral system of a
tribe living in Amazonia -- Pirah\~{a} -- has been published.
These people   use a very simple numeral system for counting: one,
two, many. For Pirah\~{a}, all quantities larger than 2 are just
`many' and such operations as 2+2 and 2+1 give the same result,
i.e., `many'. Using their weak numeral system Pirah\~{a} are not
able to see, for instance, numbers 3, 4, 5, and 6, to execute
arithmetical operations with them, and, in general, to say
anything about these numbers because in their language there are
neither words nor concepts for that. It is important to emphasize
that the records $ 2+1=\mbox{`many'}$ and $2+2= \mbox{`many'}$ are
not wrong. They are correct in their language and if one is
satisfied with the accuracy of  the answer `many', it can be used
(and \textit{is used} by Pirah\~{a}) in practice. Note that  the
result of Pirah\~{a} is not wrong, it is just \textit{inaccurate}.
Analogously, the answer `many' to the question `How many trees are
there in a park?' is correct, but its precision is low.

Thus, if one needs a more precise result than `many', it is
necessary to introduce a more powerful   numeral system   allowing
one to express the required answer in a more accurate way. By using
numeral systems where additional numerals for expressing numbers
`three' and `four' have been introduced, we can notice that within
`many' there are several objects and  numbers 3 and 4 are among
these unknown to Pirah\~{a} objects.

Our great attention to the numeral system of Pirah\~{a} is due to
the following fact: their numeral `many' gives them   such results
as
 \beq
\mbox{`many'}+ 1= \mbox{`many'},   \hspace{1cm}    \mbox{`many'} + 2
= \mbox{`many'},
    \label{Riemann29}
 \eeq
that are very familiar to us, see (\ref{Riemann30}). This
comparison shows  that we behave ourselves in front of infinity in
the same way in which Pirah\~{a} behave themselves in front of
quantities larger than 2. Thus,  our difficulty in working with
infinity is not connected to the \textit{nature of infinity
itself} but is just a result of \textit{inadequate numeral
systems} that we use to work with infinity.

In order to avoid such situations as (\ref{Riemann30}) and
(\ref{Riemann29}), a new numeral system has been proposed recently
in \cite{Sergeyev,informatica,Lagrange,chapter}. It is based on an
infinite unit of measure expressed by the numeral \G1 called
\textit{grossone}.
 A number of powerful
theoretical and applied results have been obtained  by several
authors using the new methodology. For instance,    the new
approach has been compared with  the   panorama of ideas dealing
with infinity and infinitesimals  in
\cite{Lolli,Lolli_2,MM_bijection,Sergeyev_Garro}. Then, the new
methodology has been successfully applied for studying hyperbolic
geometry (see \cite{Margenstern}), percolation (see
\cite{Iudin,DeBartolo}), fractals (see
\cite{chaos,Menger,Biology,DeBartolo}), numerical differentiation
and optimization (see \cite{DeLeone,Korea,Num_dif,Zilinskas}),
infinite series and the Riemann zeta function (see
\cite{Dif_Calculus,Riemann,Zhigljavsky}), the first Hilbert
problem and Turing machines (see
\cite{first,Sergeyev_Garro,Sergeyev_Garro_2}), cellular automata
(see \cite{DAlotto}). The usage of numerical infinitesimals opens
possibilities for creating new numerical methods having an
accuracy that is superior to existing algorithms working only with
finite numbers  (see, e.g., algorithms for solving ordinary
differential equations in \cite{ODE}).

 \begin{figure}[t]
  \begin{center}
    \epsfig{ figure = 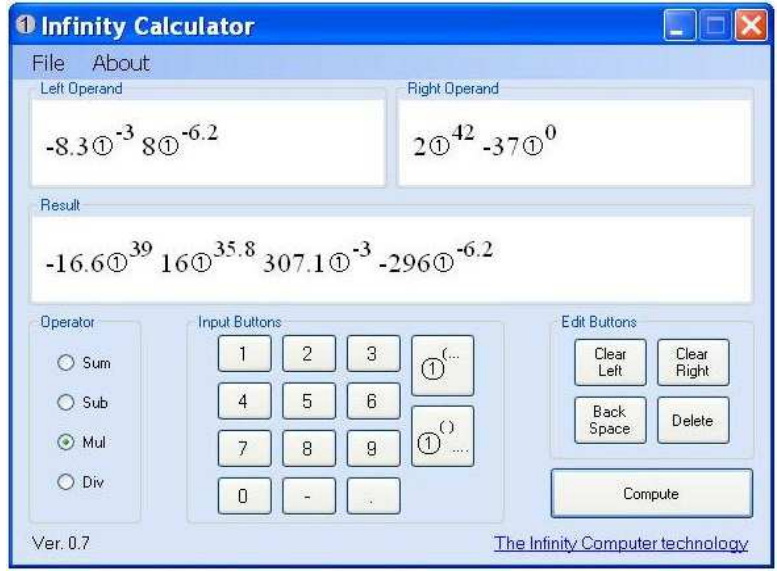, width = 4in, height = 3in,  silent = yes }
    \caption{Operation of
multiplication executed at the Infinity Calculator. The left
operand has two infinitesimal parts, the right operand has an
infinite part and a finite one; the result has two infinite and
two infinitesimal parts.}
 \label{Magic_1}
  \end{center}
\end{figure}

In particular, the Infinity Computer executing numerical
computations with infinite and infinitesimal numbers has been
patented (see \cite{Sergeyev_patent}) and its software prototype has
been constructed. This computer can be used to calculate the medal
rank $n(g_A,s_A, b_A)$ satisfying condition (\ref{medals1}) because
it works with numbers expressed in   the new positional numeral
system with the infinite base \ding{172}. A number $C$ is subdivided
into groups corresponding to powers of \ding{172}:
 \beq
  C = c_{p_{m}}
\mbox{\ding{172}}^{p_{m}} +  \ldots + c_{p_{1}}
\mbox{\ding{172}}^{p_{1}} +c_{p_{0}} \mbox{\ding{172}}^{p_{0}} +
c_{p_{-1}} \mbox{\ding{172}}^{p_{-1}}   + \ldots   + c_{p_{-k}}
 \mbox{\ding{172}}^{p_{-k}}.
\label{3.12}
       \eeq
 Then, the record
 \beq
  C = c_{p_{m}}
\mbox{\ding{172}}^{p_{m}}    \ldots   c_{p_{1}}
\mbox{\ding{172}}^{p_{1}} c_{p_{0}} \mbox{\ding{172}}^{p_{0}}
c_{p_{-1}} \mbox{\ding{172}}^{p_{-1}}     \ldots c_{p_{-k}}
 \mbox{\ding{172}}^{p_{-k}}
 \label{3.13}
       \eeq
represents  the number $C$. Numerals $c_i\neq0$, they can be
positive or negative and belong to a traditional numeral system
and are called \textit{grossdigits}. They  show how many
corresponding units $\mbox{\ding{172}}^{p_{i}}$ should be added or
subtracted in order to form the number $C$. Obviously, since all
$c_i$ are finite, it follows
 \beq
  \G1 >  c_i.
 \label{medals5}
       \eeq

Numbers $p_i$ in (\ref{3.13}) are called \textit{grosspowers}.
They are sorted in the decreasing order
\[
p_{m} >  p_{m-1}  > \ldots    > p_{1} > p_0 > p_{-1}  > \ldots
p_{-(k-1)}  >   p_{-k}
 \]
with $ p_0=0$ and, in general,  can be finite, infinite, and
infinitesimal. Hereinafter we consider only finite values of
$p_i$. Under this assumption, \textit{infinite numbers}   are
expressed by numerals having at least one $p_i>0$. They can have
several infinite parts, a finite part, and several infinitesimal
ones. \textit{Finite numbers}  are represented by numerals having
only one grosspower $ p_0=0$. In this case $C=c_0
\mbox{\ding{172}}^0=c_0$, where  $c_0$ is a conventional finite
number expressed in a traditional finite numeral system.
\textit{Infinitesimals}   are represented by numerals $C$ having
only negative grosspowers.  The simplest infinitesimal  is
$\mbox{\ding{172}}^{-1}=\frac{1}{\mbox{\ding{172}}}$ being the
inverse element with respect to multiplication for \ding{172}:
 \beq
\frac{1}{\mbox{\ding{172}}}\cdot\mbox{\ding{172}}=\mbox{\ding{172}}\cdot\frac{1}{\mbox{\ding{172}}}=1.
 \label{3.15.1}
       \eeq
Note that all infinitesimals are not equal to zero. In particular,
$\frac{1}{\mbox{\ding{172}}}>0$ because it is a result of division
of two positive numbers. Fig.~\ref{Magic_1} shows the Infinity
Calculator built using the Infinity Computer technology.

It becomes very easy  to  calculate $n(g_A,s_A, b_A)$ using records
(\ref{3.12}), (\ref{3.13}), i.e., putting \G1 instead of a finite
base $\beta$ in (\ref{medals6}). Then the number
 \beq
n(g_A,s_A, b_A)=
 g_A \G1^{2} +   s_A \G1^1 +b_A \G1^0 = g_A \G1^{2}     s_A \G1^1  b_A \G1^0 \label{medals4}
       \eeq
gives us the rank of the country  satisfying condition
(\ref{medals1}). Let us consider as an example the data
(\ref{medals3}). Since \G1 is larger than any finite number (see
(\ref{medals5})), it follows from (\ref{medals4}) that
\[
n(g_A,s_A, b_A)= 2 \cdot \G1^{2} +  0 \cdot \G1^1 + 0 \cdot \G1^0 =
2   \G1^{2}
  >
\]
\[
n(g_B,s_B, b_B)= 1 \cdot \G1^{2} +  11 \cdot \G1^1 + 0 \cdot \G1^0 =
1   \G1^{2}    11   \G1^1
\]
since
\[
2   \G1^{2} - 1   \G1^{2}    11   \G1^1 = 1   \G1^{2}    -11   \G1^1
= \G1 ( \G1 - 11) > 0.
\]
Thus, we can easily calculate the rank R1  for the data from Table~1
as follows
\[
 13 \G1^2 11  \G1^1 9 \G1^0   >
  11 \G1^2 5  \G1^1 10 \G1^0  >
    10 \G1^2 10  \G1^1 5 \G1^0  >
   9 \G1^2 7  \G1^1 12 \G1^0    >
\]
\[
    8 \G1^2 7 \G1^1     9  \G1^0 >
     8 \G1^2 6 \G1^1    5  \G1^0 >
    6 \G1^2 3   \G1^1   2 \G1^0
>  5 \G1^2  0  \G1^1    1 \G1^0
>
 \]
\[
 4 \G1^2 8  \G1^1   5 \G1^0 >   4   \G1^2  4 \G1^1 7  \G1^0.
\]

We can conclude so that the introduced calculator can be used for
computing the unofficial International Olympic Committee rank R1
numerically. Clearly, it can also be applied in all situations that
require   the lexicographic ordering not only for words with three
characters as it happens for the rank R1 but for words having any
finite number of characters, as well.

%\markboth{Bibliography}{Bibliography}
\bibliographystyle{plain}
\bibliography{XBib_Medals}
%\end{article}
\end{document}